\documentclass[12pt]{amsart}
\begin{document}

\title[ $a$-invariant and reduction numbers]
{ Bounds on the $a$-invariant and    
reduction numbers of ideals   }

\author{ Clare. D'Cruz}
\address
{Chennai Mathematical Institute \\ G. N. Chetty Road \\ T. Nagar, Chennai
600 017 India}
\email{clare@@cmi.ac.in}

\author{ Vijay Kodiyalam}
\address
{Institute of Mathematical Sciences\\  Taramani, Chennai, 600
113 India}
\email{vijay@@imsc.ernet.in}

\author{ J. K.  Verma}
\address
{Indian Institute of Technology Bombay\\ Powai, Mumbai 400 076 India}
\email{jkv@@math.iitb.ac.in}


\thanks{ AMS $2000$ classification: 13D45, 13A02.}

\thanks{Corresponding author: J. K. Verma}
\newcommand{\ncom}{\newcommand}
\ncom{\beqn}{\begin{eqnarray*}}
\ncom{\pf}{\begin{proof}}
\ncom{\eeqn}{\end{eqnarray*}}
\ncom{\beq}{\begin{eqnarray}}
\ncom{\eeq}{\end{eqnarray}}
\ncom{\been}{\begin{enumerate}}
\ncom{\eeen}{\end{enumerate}}
\ncom{\nno}{\nonumber}
\ncom{\hs}{\mbox{\hspace{.25cm}}}
\ncom{\rar}{\rightarrow}
\ncom{\lrar}{\longrightarrow}
\ncom{\Rar}{\Rightarrow}
\ncom{\Lra}{\Leftrightarrow}
\ncom{\noin}{\noindent}

\newtheorem{thm}{Theorem}[section]
\newtheorem{lemma}[thm]{Lemma}
\newtheorem{cor}[thm]{Corollary}
\newtheorem{conj}[thm]{Conjecture}
\newtheorem{pro}[thm]{Proposition}
\newtheorem{example}[thm]{Example}
\newtheorem{definition}[thm]{Definition}
\newtheorem{remark}[thm]{Remark}
\newtheorem{notation}[thm]{Notation}
\newtheorem{blank}[thm]{}

\ncom{\bt}{\begin{thm}}
\ncom{\bn}{\begin{notation}}
\ncom{\en}{\end{notation}}
\ncom{\et}{\end{thm}}
\ncom{\bl}{\begin{lemma}}
\ncom{\el}{\end{lemma}}
\ncom{\bco}{\begin{cor}}
\ncom{\eco}{\end{cor}}
\ncom{\bconj}{\begin{conj}}
\ncom{\econj}{\end{conj}}
\ncom{\ep}{\end{pro}}
\ncom{\bex}{\begin{example}}
\ncom{\eex}{\end{example}}
\ncom{\bd}{\begin{definition}}
\ncom{\ed}{\end{definition}}
\ncom{\brm}{\begin{remark}}
\ncom{\erm}{\end{remark}}
\ncom{\bb}{\begin{blank}}
\ncom{\eb}{\end{blank}}

\ncom{\sz}{\scriptsize}
\ncom{\CM}{Cohen-Macaulay }
\ncom{\sop}{system of parameters}
\ncom{\bp}{\begin{pro}}
\ncom{\eop}{\hfill{$\Box$}}
\ncom{\tfae}{the following are equivalent:}
\ncom{\mm}{minimal multiplicity }
\ncom{\f}{\frac}
\ncom{\la}{\lambda}
\ncom{\si}{\sigma}
\ncom{\ssize}{\scriptsize}
\ncom{\al}{\alpha}
\ncom{\be}{\beta}
\ncom{\Si}{\Sigma}
\ncom{\kbar}{\overline{\kappa}}
\ncom{\bib}{\bibitem}
\ncom{\sst}{\subset}
\ncom{\sms}{\setminus}
\ncom{\seq}{\subseteq}
\ncom{\est}{\emptyset}
\ncom{\bighs}{\hspace{.5 cm}}
\ncom{\ulin}{\underline}
\ncom{\ol}{\overline}
\ncom{\sta}{\stackrel}
\ncom{\bop}{\bigoplus}

\ncom{\depth}{\mbox{depth}}
\font\Bbb=msbm10 scaled 1200
\font\Bb=msbm10 scaled 800
\font\Bbbb=msbm10 scaled 1800
\newcommand{\RR}{\mbox{\Bbb R}}
\newcommand{\CC}{\mbox{\Bbb C}}
\newcommand{\C}{\mbox{\Bb C}}
\newcommand{\GG}{\mbox{\Bbb G}}
\newcommand{\ZZ}{\mbox{\Bbb Z}}
\newcommand{\QQ}{\mbox{\Bbb Q}}
\newcommand{\NN}{\mbox{\Bbb N}}
\newcommand{\FF}{\mbox{\Bbb F}}
\newcommand{\KK}{\mbox{\Bbb K}}
\newcommand{\HH}{\mbox{\Bbb H}}
\newcommand{\ga}{\mbox{\german a}}
\newcommand{\gb}{\mbox{\german b}}

\def\m{{\frak m}}
\def\M{{\frak M}}

\begin{abstract}
Let $R$ be a $d$-dimensional standard graded ring over an
Artinian local ring. Let $\M$ be the unique
maximal homogeneous ideal of $R$. Let $h^i(R)_n$ denote the length of the
$nth$ graded component of the local cohomology module $H^i_{\M}(R).$    
Define the Eisenbud-Goto invariant 
$EG(R)$ of $R$ to be the number 
${\sum_{q=0}^{d-1}{d-1 \choose q}h^q_{\M}(R)_{1-q}}.$
We prove that
the $a$-invariant $a(R)$ of the top local cohomology module 
$H_{\M}^d(R)$ satisfies the inequality: 
$a(R) \leq e(R)-\ell(R_1)+(d-1)(\ell(R_0)-1)+EG(R).$ This
 bound is used to get  upper bounds  for the reduction number of
an $\m$-primary ideal $I$ of a \CM local ring $(R,\m  )$, 
when the associated graded ring of $I$ has depth at least $d-1.$ 
\end{abstract}  

\maketitle

\section{Introduction}

 Let $R=\bigoplus_{n=0}^{n=\infty} R_n$ be a $d$-dimensional
standard  graded ring over an Artinian local  ring $R_0$.
 Let $\M$ be the maximal homogeneous ideal
of $R$.  Let  $H^i_{\M}(R)$ denote the $i$-th local cohomology module of $R$
 with
respect to  $\M$. For a graded module $M$, we use $[M]_n$ or $M_n$ to denote the 
$nth$ graded component of $M.$  The    $a$-invariant of $R,$  introduced 
in \cite{gotow}, is defined as  
$$a(R)={\max}\{n \left |[H^d_{\M}(R)]_n  \neq 0 \right. \}.$$
The objective  of this paper is to give a bound
for the $a$-invariant of $R$ in terms of lengths of graded 
components of local cohomology modules and use it to get 
bounds for reduction numbers of ideals. Let $\ell(M)$ denote length
of a module $M.$
We set $\ell\left(\left[H^q_{\M}(R)\right]_{1-q}\right)=h^q(R)_{1-q}$
for all $q \geq 0.$
 To state our bound for the 
$a$-invariant we define the Eisenbud-Goto invariant  $EG(R)$ of  $R$ to be 
the number $$  
EG(R)={\sum_{q=0}^{d-1}{d-1 \choose q}h^q (R)_{1-q}}.
$$
The main result of the paper is: 
\bt
Let $R=\bigoplus^{\infty}_{n=0}R_n$ be a $d$-dimensional standard 
graded algebra over an artinian local ring $R_0$  
with multiplicity $e(R).$ Then
$$
a(R) \leq e(R)-\ell(R_1)+(d-1)(\ell(R_0)-1)+EG(R).
$$
\et
\noindent
 Eisenbud and Goto \cite{eg} showed that if $R_0$ is a field then
$$
e(R) \geq 1+\mbox{codim}(R)-EG(R).
$$
They showed that if equality holds in the above inequality
then $R/H^0_{\M}(R)$ has linear resolution. 
To state our bounds for reduction numbers we recall some basic concepts
about reductions of ideals.
Let $(R,\m)$ be a local ring. Let $J \subset I$ be ideals of $R$. The
ideal $J$ is called a reduction of $I$ if there exists  an $n \in \NN$
such that $JI^n=I^{n + 1}$  \cite{nr}. Among the reductions of $I,$
the smallest ones with respect to inclusion are called minimal
reductions of $I$. If $R/\m  $ is infinite then any minimal reduction
of $I$ is minimally generated by as many elements as the Krull
dimension of the fiber cone  $F(I):= \bop_{n = 0}^{n =  \infty} I^n/\m
I^n$.  The {\it reduction number,} $r_J(I),$ of $I$ with respect to a
minimal  reduction $J$ is the least integer $n$ for which $JI^n = I^{n
  + 1}$. When $R/\m  $ is infinite, the reduction number of $I$ is
defined as the minimum of the reduction numbers $r_J(I)$ where $J$
varies over all the minimal reductions of $I.$
Let $G(I) := \oplus_{n \geq   0} {I}^{n} / {I}^{n+1}$ be the
associated  graded ring  of an  ideal $I$.  
 Let $\gamma(I)$ denote the depth of the irrelevant ideal $G_+$
of  $G(I)$. If $(R,\m )$ is a Cohen-Macaulay local
 ring, $I$ is an $\m$-primary ideal and  $\gamma(I) \geq d-1$, 
then $r(I) = a(G(I))  + d$  \cite{m2}. The {\em Ratliff-Rush closure} of an 
ideal $I,$ $\widetilde{I},$ is the stable value of the sequence of the ideals
$\{I^{n+1}:I^n\}.$ We will obtain the following
bounds for $r(I)$ as an application of the main theorem:

\bt Let $(R,\m)$ be a $d$-dimensional \CM local ring 
with infinite residue field. Let
$I$ be an $\m$-primary ideal with $\gamma(I) \geq d-1.$ 
Let $J$ be any minimal reduction of $I.$  
\begin{enumerate}
\item Let $d=1.$ Then 
$r(I) \leq e(I)-
      \left(\ell(I/(I\cap \widetilde{I^2})-1\right)\leq e(I).$
 \item Let $d=2.$ Put $X= \mbox{Proj}(G(I)).$ Then 
           $r(I)\leq 1+e(I)-\ell(I/I^2)+\ell(H^0(X,{\mathcal O}_X)).$
 \item Let $d \geq 3.$ Then  
 $
 r(I) \leq 1+\ell(I^2/JI)+ h^{d-1}(G(I))_{2-d}.
 $
\end{enumerate}
\et
We will show by an example that our bounds for the 
$a$-invariant and reduction number are sharp.


\section{A bound on the  a-invariant of  standard graded algebras}
\label{regularity}

In this section we prove our bound on the  $a$-invariant
 of a standard  graded algebra  $R$ over an Artinian local ring
 $R_0$.

\bt
\label{main}
Let $R =\bigoplus_{n=0}^{n = \infty} R_n$ be a $d$-dimensional
standard graded algebra  over an Artinian local ring $R_0.$ 
Then
\beq
\label{onee}
          a(R) 
& \leq &  e(R) - \ell(R_1) + (d-1) (\ell(R_0) -1)
               + EG(R).
\eeq
\et
\pf
We may assume without loss of generality that the 
residue field of $R_0$ is infinite. We  prove the theorem by
induction on $d$.  Let $d = 0.$ Then
$$  
e(R)=\ell(R_0)+\ell(R_1)+\cdots+\ell(R_m),
$$
where $m=a_0(R).$ Thus
$$e(R)-\ell(R_1)-\ell(R_0)+1=1+\ell(R_2)+\cdots+\ell(R_m) \geq m.$$
Let $R$ be \CM and pick a degree one nonzerodivisor $x$ to see that 
\beqn
a(R) &=& a(R/xR)-1 \\ 
     &\leq & e(R/xR)-\ell\left(\left[R/xR\right]_1\right)+
      (d-2)(\ell(R_0)-1)-1 \\ 
     & = & e(R)-\ell(R_1)+\ell(R_0)+(d-2)(\ell(R_0)-1)-1 \\
     & = & e(R)-\ell(R_1)+(d-1)(\ell(R_0)-1). 
\eeqn
Now let $d =  1$. If $R$ is \CM, we are done by the above argument. 
So let  depth$(R)=0.$ Then $S:=R/H^0_\M(R)$ is \CM, $e(S)=e(R)$
and $a(R)=a(S).$ Hence  
$$ 
a(R)=a (S) \leq e(S) - \ell(S_1)=e(R)-\ell(R_1)+h^0(R)_1 
.$$

Suppose  $d \geq  2$. Let $x \in R_1 $ be a superficial
element. We first  prove that for a degree one superficial 
element in $R,$
$$EG(R/xR) \leq EG(R).$$
Since $x$ is superficial of degree one,
\beqn
H_{\M}^{i}(0:_R x) =   \left\{ \begin{array}{ll}
                        (0 :_R x) & \mbox{if} \; i = 0 \\
                         0 & \mbox{if otherwise.} \\
                                           \end{array}
                                           \right. .
\eeqn
Hence from the short exact sequence 
\beqn
       0 
\lrar (0:_R x)
\lrar R 
\lrar \f{R}{(0:_R x)} 
\lrar 0
\eeqn
we get 
$H_{\M}^{i}(R/(0:_R x)) = H_{\M}^{i}(R)
$
for all $i \geq 1$. 
From the exact sequence 
\beqn
        0 
\lrar  \f{R}{(0:_R x)}(-1) 
\lrar R 
\lrar  \f{R}{xR} \lrar 0
\eeqn
we get the long exact sequence
\beqn
  \cdots \lrar              [H_{\M}^{i}(R)]_n 
  \lrar   [H_{\M}^{i}(R/xR)]_n
  \lrar   [H_{\M}^{i+1}(R)]_{n-1} \\
  \lrar [H_{\M}^{i+1}(R)]_n
  \lrar \cdots .
\eeqn
Hence for all $i \geq 0$,
$$
       h^{i}(R/xR)_{n} 
      \leq h^{i}(R)_{n} 
       +  h^{i+1}(R)_{n-1} .
$$
Hence
\beqn
EG(R/xR)&=& \sum_{q=0}^{d-2}{d-2 \choose q}h^q(R/xR)_{1-q}\\
       &\leq & \sum_{q=0}^{d-2}{d-2 \choose q}
        \left[h^q(R)_{1-q}+h^{q+1}(R)_{-q}\right]\\
       &= & \sum_{q=0}^{d-2}{d-2 \choose q}h^q(R)_{1-q}+
                \sum_{q=1}^{d-1}{d-2 \choose q-1} h^{q}(R)_{1-q}\\
       &=& \sum_{q=0}^{d-1}{d-1 \choose q} h^{q}(R)_{1-q}\\
       &=& EG(R).
\eeqn
Therefore
\beqn
a(R) &\leq & a(R/xR)-1 \;\;\; \mbox{by \cite{trung}}\\
     &\leq&e(R/xR)-\ell(R/xR)_1+(d-2)
           \left(\ell\left(\left[R/xR\right]_0\right)-1\right)+EG(R)-1 \\
     &=&   e(R)-\ell(R_1)+\ell(R_0)-\ell\left((0:x)_{R_0}\right)+
           (d-2)\left(\ell(R_0)-1\right) \\
     & + & EG(R)-1 \\
     &\leq& e(R)-\ell(R_1)+
           (d-1)\left(\ell(R_0)-1\right)+EG(R). \
\eeqn
\end{proof}

\noindent
We now demonstrate that the bound in  Theorem~\ref{main} is sharp. 

\bex\label{counter}
{\em
Let $k$ be a field and $x,y,a,b,c,d$ be
indeterminates. Consider the ideal $I = (x^3, x^2y^4, xy^5, y^7)$ in the
polynomial ring $S = k[x,y].$ 
Using Hilbert series we show that 
$F(I) \simeq k[a, b, c, d]/(bd, bc, b^2, c^3).$ Consider the ring
 homomorphism $\phi : R = k[a, b, c, d] \lrar F(I)$ defined by
$$
\phi(a)  =  \overline{x^3},\;\;  \phi(b)  =  \overline{x^2y^4},\;\;
\phi(c)  =  \overline{xy^5},\;\; \mbox{and}\;\; \phi(d) = \overline{y^7}.
$$ 
Here the overbar indicates the image in the first graded component
of $F(I).$ Let $L=\mbox{ker}~\phi$. The equations
\begin{eqnarray*}
  (x^2y^4)(y^7)    & =&  (xy^5)^2y, \\ 
  (x^2y^4)(xy^5) & = & (x^3y^7)y^2, \\
  (x^2y^4)^2     & = & (x^3)(xy^5)y^3, \\ 
  (xy^5)^3       & = & x^3(y^7)^2y,
\end{eqnarray*}
show that $N = (bd, bc, b^2, c^3) \subset L.$ To show that $N = L$, we show
that $R/N$ and $R/L$ have same Hilbert series. We denote the Hilbert series
of a graded algebra $G$ by $H(G, \lambda).$ By the propositions 2.3 and
2.6 of  \cite{H} we find that $\mu(I^n) = 3n + 1$ for all $n \geq 0.$ Here
$\mu$ denotes the minimum number of generators. This shows that
$H(F(I),\lambda) = (1 + 2\lambda)/(1 - \lambda)^2.$ 
By the well known "divide and
conquer strategy" for finding Hilbert series of quotients of polynomial rings by monomial ideals we get,
$H(R/N,\lambda) =  H(F(I),\lambda).$ Thus $F(I) \cong R/N.$ Therefore $F(I)$ is a
two - dimensional ring with depth one. Notice that 
$N = (b, c^3) \cap (c, d, b^2).$ Put $J = (b, c^3)$ and $K = (c, d, b^2).$
In order to get the desired information about local cohomology of $F(I)$,
consider  the exact sequence :
$$
     0 
\lrar F(I) 
\lrar R/J \bop R/K 
\lrar R/(J + K) \lrar 0.
$$ 
Hence we get  the following  long exact sequence of local cohomology 
modules with respect to the maximal homogeneous ideal $\m   = (a, b, c, d) :$

$ 
      0 
\lrar H_\m  ^1(F(I)) 
\lrar H_\m  ^1(R/K) $ 

$\lrar H_\m  ^1(R/(J + K)) 
\lrar  H_\m  ^2(F(I)) 
\lrar H_\m  ^2(R/J) 
\lrar 0.
$

\noindent
We now show that $a(F(I))=0$ and $h^1(F(I))_0=1.$ Since
$$R/J \simeq k[a,c,d]/(c^3),\;\; R/(J+K) \simeq k[a] 
\;\;\mbox{and} \;\;R/K \simeq k[a,b]/(b^2),$$
by using that fact that $a(R/(f))=a(R)+\mbox{deg}(f)$ 
for a homogeneous regular 
element $f$ of a graded algebra $R$, we conclude that 
$a(R/J)=0,\;\;a(R/(J+K))=-1 
\;\;\mbox{and} \;\;a(R/K)=0.$ Thus $a(F(I))=0.$ By \cite[Theorem 
4.4.3]{bh}, we 
get 
$$h^1(F(I))_0=h^1(R/K)_0=P(R/K,0)-H(R/K,0)=2-1=1.$$  
Substituting these values in  
(\ref{onee}) we observe that equality holds. Therefore
(\ref{onee}) is sharp. }
\qed
\eex

\section{Bounds on reduction numbers}
\label{graded}
In this section we will use the bound on the $a$-invariant obtained 
in the previous section
to provide bounds on reduction numbers. By \cite{trung} and  
\cite{m2}, we know that $r(I)=a(G(I))+d$ where $(R,\m)$ is a 
\CM local ring of dimension $d$ and $\gamma(I) \geq d-1.$ We will consider the cases where 
$d=1,\;\;d=2\;\;$ and $d \geq 3$ separately.
In the next result we will need the formula:
$[H^0_{G_{ + }}(G(I))]_n = (I^n\cap \widetilde{I^{n + 1}})/I^{n + 1}$ for
all $n \geq 0$  \cite{hjls}.

\bp
\label{associated1}
Let $(R,\m  )$ be a one-dimensional \CM local ring. Let $I$ be an
$\m$-primary ideal. Then 
$$ r(I) \leq e(I) - \left[\ell (I/(I\cap \widetilde{I^2}) -1\right]
 \leq e(I).$$ 

\ep

\pf 
Since $d=1$, 
\beqn
 a(G(I)) &\leq &   e(I)- \ell(I/I^2) + h^0(G)_1 \\
     &=    &   e(I)- \ell(I/I^2) + \ell((I\cap \widetilde{I^2})/I^2) \\
     &=    & e(I)-  \ell(I/(I\cap \widetilde{I^2})).
\eeqn
Hence $r(I)=a(G(I))+1 \leq 1+ e(I)-  \ell(I/(I\cap \widetilde{I^2}).$
If $ \ell(I/(I\cap \widetilde{I^2})=0$ then 
$I \subseteq \widetilde{I^2}$. But $\widetilde{I^2}= I^{n+2}:I^n$
 for large $n.$ Hence $I^{n+1}=I^{n+2}.$ This is a contradiction. Hence
 $\ell(I/(I\cap \widetilde{I^2})) \geq 1.$ Thus we obtain the classical 
bound $r(I) \leq e(I).$

\end{proof}

\bex
{\em Let $k$ be a field and $t$ be an indeterminate. Put 
$R = k[[t^4, t^5, t^6, t^7]]$ and 
$I = (t^4, t^5, t^6).$  Let $\m  $ denote the unique
maximal ideal of $R.$ Let 
$G$ denote the associated graded ring $G(I)$ 
of $I$. Then $G$ is not \CM since $t^7 I\subset I^2.$   To find
the  associated Ratliff-Rush ideal of $I$ notice that $I^2 = \m  ^2.$
Since $r(\m)=1,$ 
the associated graded ring $G(\m)$ is \CM by
\cite{sally}. Therefore all powers of $\m  $ are
Ratliff-Rush. Hence, $\widetilde{I^2} = \widetilde{\m^2} = \m^2 =
I^2.$  Hence $ (I\cap \widetilde{I^2}/I^2) = 0.$ Therefore $r(I) \leq
1 +  e(I) - \ell(I/I^2)  =  2.$ It can be checked that 
$r_{(t^4)}(I)=2.$ Therefore the bound in the above result is sharp.}

\eex

\bp
\label{associated}
Let $I$ be an $\m  $-primary ideal of a  two dimensional \CM local ring with
$\gamma(I)\geq 1.$ Let $X={\mbox Proj}\;G(I).$  Then 
$$r(I) \leq 1+e(I)-\ell(I/I^2)+ \ell(H^0(X,{\mathcal O}_X)).$$
\ep

\pf 
Since 
$\gamma(I) \geq 1$, 
$$
     r(I)
\leq 1+ e(I)  -  \ell(I/I^2) + \ell(R/I) 
+     h^{1}(G)_{0}.
$$  
By the exact sequence
$$ 0\longrightarrow H^0_{G_+}(G) \longrightarrow G
    \longrightarrow \bigoplus_{n \in \ZZ}H^0(X,{\mathcal O}_X(n))
    \longrightarrow H^1_{G_+}(G)\longrightarrow 0
$$
we get, by taking the $0th$ component of all the modules in 
the above exact sequence:
$$
 \ell(H^0(X,{\mathcal O}_X))-\ell(R/I)= h^1(G)_0
$$
Putting this in the above bound for $r(I)$ we get the desired upper bound.
\end{proof}

\bp 
Let $(R,\m)$ be a \CM local ring of dimension $d \geq 3$ with $R/\m$
infinite. Let     $I$ be an 
$\m$-primary ideal with a minimal reduction $J$ and $\gamma(I) \geq d-1.$ Then
$$ r(I) \leq 1+\ell(I^2/JI) + h^{d-1}(G)_{2-d}.$$
\ep
\pf 
Since $\gamma(I) \geq d-1,$ by the Theorem \ref{main} 
\beqn
r(I)& \leq & e(I)-\ell(I/I^2)+(d-1)(\ell(R/I)-1)+ h^{d-1}(G)_{2-d}+d \\
    & =    & \ell(R/J)-\ell(R/I^2)+
             \ell(R/I)+(d-1)\ell(R/I) \\
    &- &(d-1)+d+ h^{d-1}(G)_{2-d}\\
    & =    &  \ell(R/J)-\ell(R/I^2)+\ell(J/JI)+1+ h^{d-1}(G)_{2-d}\\
    &  =   &  1+\ell(I^2/JI)+ h^{d-1}(G)_{2-d}\\
\eeqn
\end{proof}

{\bf Acknowledgements:}
The first author wishes to acknowledge the hospitality of the  Indian
Institute of Technology Bombay, where part of the paper was done. 
The last named author wishes to thank the  Institute
of Mathematical Sciences where this work was initiated.

\end{document}